\documentclass[12pt]{article}
\usepackage{mathrsfs}
\usepackage{amssymb} \textwidth 150mm \textheight 220mm \oddsidemargin
15pt \evensidemargin 0pt \topmargin 0cm \headsep 0.3cm

\usepackage{amsmath}
\usepackage{graphicx}
\usepackage{indentfirst}
\usepackage{cite}
\usepackage{float}
\usepackage{dsfont}
\usepackage{CJK}
\usepackage{multirow}
\usepackage{makecell}

\usepackage{subfigure}

\begin{document}
\title{\Large\bf{Limit cycles appearing from the perturbation of a cubic isochronous center }}
\author{{Jihua Yang\thanks{Corresponding author.\newline\ E-mail addresses: jihua1113@163.com(J. Yang), qipeng\_zhang@outlook.com(Q. Zhang) },\ \ Qipeng Zhang}
\\ {\small \it School of Mathematical Sciences, Tianjin Normal University, Tianjin 300387, China}\\
 }
\date{}
\maketitle \baselineskip=0.9\normalbaselineskip \vspace{-3pt}
\noindent
{\bf Abstract}\, For a polynomial differential system
$$\dot{x}=-y+\sum\limits_{i+j=3}\alpha_{i,j}x^iy^j,\quad \dot{y}=x+\sum\limits_{i+j=3}\beta_{i,j}x^iy^j,$$
 Pleshkan (Differ. Equations, 1969) proved that the origin is an isochronous center of this system iff it can be brought to one of $S^*_1$,  $S^*_2$, $S^*_3$ or $S^*_4$. The bifurcation of limit cycles for these four types of isochronous differential systems have not yet been studied, except for $S^*_1$. This paper is devoted to study the limit cycle problem of $S^*_2$ when we perturb it with an arbitrary polynomial vector field. An upper bound of the number of limit cycles is obtained using the Abelian integral.
\vskip 0.2 true cm
\noindent
 \section{Introduction and main result}
 \setcounter{equation}{0}
\renewcommand\theequation{1.\arabic{equation}}

The number of limit cycles and their relative positions play a vital role in determining the global phase portrait of a differential system. Moreover, they possess significant practical significance as they reflect the widespread periodic oscillation phenomena in the real world. Isochronism appears in numerous practical problems and manifests itself as synchronous oscillations in the field of physics. Therefore, studying the bifurcation of limit cycles of differential systems with isochronous centers is of great practical application value \cite{C,CL06,CL07,CL08,CG99}. Another important reason why the bifurcation of limit cycles of differential systems has drawn special attention is its close connection with the renowned Hilbert's 16th problem, which is to study the maximum number of limit cycles generated by two-dimensional differential systems, including isochronous differential systems, and their relative positions \cite{H}.

It is well known that the problem of isochronicity appears only for the non-degenerate centers. A polynomial differential system of degree $n$ with a non-degenerate center at the origin can be taken the form
\begin{eqnarray}
\dot{x}=-y+\sum\limits_{i+j=2}^n\alpha_{i,j}x^iy^j,\ \dot{y}=x+\sum\limits_{i+j=2}^n\beta_{i,j}x^iy^j,\ a_{i,j},b_{i,j}\in\mathds{R},
\end{eqnarray}
in an appropriate coordinate system and upon a time rescaling. When $n=2$, Loud concluded \cite{L} that the origin is an isochronous center for system (1.1) if and only if, after a linear transformation and a time rescaling, it can be transformed into  one of the following systems
 \begin{eqnarray*}
\begin{aligned}
&S_1:\ \dot{x}=-y+x^2-y^2,\ \ \ \dot{y}=x+2xy,\\
&S_2:\ \dot{x}=-y+x^2,\ \ \ \dot{y}=x+xy,\\
&S_3:\ \dot{x}=-y-\frac{4}{3}x^2,\ \ \ \dot{y}=x-\frac{16}{3}xy,\\
&S_4:\ \dot{x}=-y+\frac{16}{3}x^2-\frac{4}{3}y^2,\ \ \ \dot{y}=x+\frac{8}{3}xy.\\
\end{aligned}
\end{eqnarray*}
Regarding the bifurcation of limit cycles in the aforementioned four classes of quadratic isochronous differential systems, extensive literature has conducted in-depth research and achieved relatively comprehensive results. Under quadratic polynomial perturbations, Chicone and Jacobs \cite{CJ} proved that at most one limit cycle can bifurcate from each period annulus of $S_1$, while at most two limit cycles can bifurcate from each period annulus of $S_2$, $S_3$ and $S_4$. Under $n$th-degree polynomial perturbations, Li et al. \cite{LLLZ} applied the  Abelian integral combined with analytical techniques to derive upper bounds on the number of limit cycles for these systems. Cen et al. \cite{CLYZ} investigated the limit cycle bifurcations of them under the non-smooth perturbation of $n$th-degree polynomials. Yang \cite{Y} utilized the Melnikov function in conjunction with the Picard-Fuchs equation to obtain upper bounds on the number of limit cycles of $S_1$ and $S_2$ under $n$th-degree polynomial perturbations with two switching lines.

For the linear plus cubic homogeneous polynomial differential system of the form
\begin{eqnarray}
\dot{x}=-y+\sum\limits_{i+j=3}\alpha_{i,j}x^iy^j,\quad \dot{y}=x+\sum\limits_{i+j=3}\beta_{i,j}x^iy^j,
\end{eqnarray}
Pleshkan \cite{P} attested that the origin is an isochronous center if and only if this system can be brought to one of the following systems:\begin{eqnarray*}
\begin{aligned}
&S^*_1:\ \dot{x}=-y+x^3-3xy^2,\ \ \ \dot{y}=x+3x^2y-y^3,\\
&S^*_2:\ \dot{x}=-y+x^3-xy^2,\ \ \ \dot{y}=x+x^2y-y^3,\\
&S^*_3:\ \dot{x}=-y+3x^2y,\ \ \ \dot{y}=x-2x^3+9xy^2,\\
&S^*_4:\ \dot{x}=-y-3x^2y,\ \ \ \dot{y}=x+2x^3-9xy^2.\\
\end{aligned}
\end{eqnarray*}
 Grau and Villadelprat \cite{GV} proved that at most two critical periods bifurcate from each period annulus of  $S^*_1$, $S^*_2$, $S^*_3$ and $S^*_4$ when they are perturbed by polynomials of degree 3, and this bound can be reached. The limit cycle problem of $S^*_1$ under perturbations of polynomials of degree $n$ was considered by Gasull et al. in \cite{GLLZ}. As far as we know, there are no more results regarding the number of limit cycles of the aforementioned isochronous systems, because the first integrals of them are quite involved.

In the present paper, we intend to investigate the limit cycle problem of the system $S_2^*$ under perturbation of $n$th-degree polynomials in $x$ and $y$. That is,
\begin{eqnarray}
\dot{x}=-y+x^3-xy^2+\varepsilon f(x,y),\ \ \dot{y}=x+x^2y-y^3+\varepsilon g(x,y),
\end{eqnarray}
where $0<|\varepsilon|\ll1$, $$f(x,y)=\sum\limits_{i+j=0}^na_{i,j}x^iy^j,\ g(x,y)=\sum\limits_{i+j=0}^nb_{i,j}x^iy^j,\ a_{i,j},b_{i,j}\in\mathbb{R},i,j\in\mathbb{N}.$$
 System $S_2^*$ has a first integral
\begin{eqnarray}
H(x,y)=(x^2+y^2)(1+2xy)^{-1}=h,\ h\in(0,1),
\end{eqnarray}
with integrating factor $\mu(x,y)=2(1+2xy)^{-2}$. The symmetric period annulus, surrounding the isochronous centre at $(x,y)=(0,0)$ (corresponding to $h=0$), is denoted by $\Gamma_h=\{(x,y): H(x,y)=h,\, h\in(0,1)\}$.

It is evident that the total number of zeros of the Abelian integral
\begin{eqnarray}I(h)=\oint_{\Gamma_h}\mu(x,y)g(x,y)dx-\mu(x,y)f(x,y)dy\end{eqnarray}
provides an upper bound for the number of limit cycles of system (1.3) bifurcating from the corresponding period annulus, and the existence of multiple simple zeros provides a lower bound of the number of limit cycles. Our main result is the following theorem.
\vskip 0.2 true cm

\noindent
{\bf Theorem 1.1}\, {\it  If the Abelian integral (1.5) is not identically zero, then, for any $n\geq1$, system (1.3) has at most $4[\frac{n+1}{2}]+1$ limit cycles bifurcating from the period annulus. }

\section{Proof of Theorem 1.1}
 \setcounter{equation}{0}
\renewcommand\theequation{2.\arabic{equation}}

The first goal of this section is to obtain the algebraic structure of the Abelian integral $I(h)$. A straightforward calculation using (1.4) gives rise to
\begin{eqnarray}\begin{aligned}
I(h)=2h^2\Big(\sum\limits_{i+j=0}^nb_{i,j}\oint_{\Gamma_h}\frac{x^iy^j}{(x^2+y^2)^2}dx-\sum\limits_{i+j=0}^na_{i,j}\oint_{\Gamma_h}\frac{x^iy^j}{(x^2+y^2)^2}dy\Big).
\end{aligned}\end{eqnarray}
Introducing the change of variables $x=r\cos\theta,y=r\sin\theta$, one has, in light of (1.4)
\begin{eqnarray}
r(\theta)=\sqrt{\frac{h}{1-h\sin(2\theta)}},
\end{eqnarray}
For a positive number $\delta$  small enough, applying the Green's formula to (2.1), one has
\begin{eqnarray}\begin{aligned}
I(h)=&h^2\sum\limits_{i+j=1}^{n+1}c_{i,j}\iint_{D_{h,\delta}}\frac{x^iy^j}{(x^2+y^2)^3}dxdy\\
&-2h^2\Big(\sum\limits_{i+j=0}^nb_{i,j}\oint_{r=\delta}\frac{x^iy^j}{(x^2+y^2)^2}dx-\sum\limits_{i+j=0}^na_{i,j}\oint_{r=\delta}\frac{x^iy^j}{(x^2+y^2)^2}dy\Big)\\
=&h^2\sum\limits_{i+j=1}^{n+1}c_{i,j}\int_0^{2\pi}\int_\delta^{r(\theta)}r^{i+j-5}\cos^i\theta\sin^j\theta dr d\theta\\
&-2h^2\Big(\sum\limits_{i+j=0}^nb_{i,j}\oint_{r=\delta}\frac{x^iy^j}{(x^2+y^2)^2}dx-\sum\limits_{i+j=0}^na_{i,j}\oint_{r=\delta}\frac{x^iy^j}{(x^2+y^2)^2}dy\Big)\\
=&h^2\sum_{\substack{i+j=1,\\ i+j\neq4}}^{n+1}c_{i,j}\int_0^{2\pi}\frac{r(\theta)^{i+j-4}}{i+j-4}\cos^i\theta\sin^j\theta d\theta\\
&+h^2\sum\limits_{i+j=4}c_{i,j}\int_0^{2\pi}\ln r(\theta)\cos^i\theta\sin^j\theta d\theta-h^2C_\delta,\\
\end{aligned}\end{eqnarray}
where $c_{i,j}$ are constants and can be chosen arbitrarily, and
$$\begin{aligned}C_\delta=&2\sum\limits_{i+j=0}^nb_{i,j}\oint_{r=\delta}\frac{x^iy^j}{(x^2+y^2)^2}dx-2\sum\limits_{i+j=0}^na_{i,j}\oint_{r=\delta}\frac{x^iy^j}{(x^2+y^2)^2}dy\\
&+\sum_{\substack{i+j=1,\\ i+j\neq4}}^{n+1}c_{i,j}\int_0^{2\pi}\frac{\delta^{i+j-4}}{i+j-4}\cos^i\theta\sin^j\theta d\theta+\sum\limits_{i+j=4}c_{i,j}\int_0^{2\pi}\ln\delta\cos^i\theta\sin^j\theta d\theta.\end{aligned}$$
It is clear that $C_\delta$ does not depend on $h$. Hence one can eliminate it by differentiating with respect to $h$.

 From the above analysis, one can obtain
\begin{eqnarray}
I(h)=h^2\big(\sum\limits_{i+j=1}^{n+1}c_{i,j}I_{i,j}(h)-C_\delta\big)\triangleq h^2I_1(h),
\end{eqnarray}
where $$I_{i,j}(h)=\begin{cases}
\frac{1}{i+j-4}\int_0^{2\pi}r(\theta)^{i+j-4}\cos^i\theta\sin^j\theta d\theta,\ i+j\neq4,\\
\int_0^{2\pi}\ln r(\theta)\cos^i\theta\sin^j\theta d\theta,\qquad\qquad\ i+j=4.
\end{cases}$$
Since $I(h)$ and $I_1(h)$ have the same number of zeros on (0,1), it suffices to obtain the expression of $I_1(h)$.  To this end, we state a result of Gasull, Li, Llibre and Zhang published in \cite{GLLZ} below which we will use to compute $I_{i,j}(h)$ in (2.4).
\vskip 0.2 true cm

\noindent
{\bf Lemma 2.1}\, {\it Let $f(x)$ be a continuous function and let $i,j\geq0$ be integers. Then the following conclusions hold:}

\noindent
(i) {\it If $i+j$ is odd, then $$\int_0^{2\pi} f(\sin2\theta)\cos^i\theta\sin^j\theta d\theta=0.$$}

\noindent
(ii) {\it If $i+j=2N$ even, then there exist real constants $d_0,d_1,\cdots,d_N$, such that
 \begin{eqnarray*}\begin{aligned}
\int_0^{2\pi}f(\sin2\theta)\cos^i\theta\sin^j\theta d\theta=\sum\limits_{k=0}^N d_k\int_0^{2\pi} f(\sin\theta)\sin^k\theta d\theta.
\end{aligned}\end{eqnarray*}}

Thanks to Lemma 2.1, one gets that $I_{i,j}(h)=0$ for $i+j$ odd, and for $i+j=2N$ even
\begin{eqnarray}\begin{aligned}
I_{i,j}(h)=\begin{cases}
\sum\limits_{k=0}^Nd_kh^{N-2}\int_{-\pi}^\pi\frac{\sin^k\theta}{(1-h\sin\theta)^{N-2}}d\theta,\ i+j=2N\neq4,\\
\sum\limits_{k=0}^2d_k\int_{-\pi}^\pi\sin^k\theta\ln\frac{h}{1-h\sin\theta}d\theta,\quad\ i+j=4.
\end{cases}\end{aligned}\end{eqnarray}
Therefore, $$I_1(h)=\sum\limits_{N=1}^{[\frac{n+1}{2}]}\sum\limits_{i+j=2N}c_{i,j}I_{i,j}(h)-C_\delta.$$
For abbreviation we use the notation
$$J_N(h)=\sum\limits_{i+j=2N}c_{i,j}I_{i,j}(h).$$

Our task now is to calculate the integrals $I_{i,j}(h)$ that appear in $J_N(h)$. When $N=1$, according to (2.5), one has
\begin{eqnarray}\begin{aligned}
I_{i,j}(h)=\sum\limits_{k=0}^1d_k\int^\pi_{-\pi}\sin^k\theta \frac{1-h\sin\theta}{h}d\theta=2\pi d_0 h^{-1}-\pi d_1.
\end{aligned}\end{eqnarray}

When $N=2$, in view of (2.5), one obtains
 $$I_{i,j}(h)=\sum\limits_{k=0}^2d_k \Big[\ln h \int^\pi_{-\pi}\sin^k\theta d\theta-A_k(h)\Big],$$
 where $$A_k(h)=\int^\pi_{-\pi}\sin^k\theta \ln(1-h\sin\theta)d\theta,\ k=0,1,2.$$
We are now turning to deal with the integrals $A_k(h).$ However, it is challenging to calculate them directly, so we compute their derivatives instead. To this effect, we denote
\begin{eqnarray}
L_k(h)=\int_0^{2\pi}\frac{1}{(1-h\sin\theta)^k}d\theta,\ h\in(0,1),k\in\mathbb{N}.
\end{eqnarray}
This integral was studied by Gasull, Li and Torregrosa in \cite{GLT}, and they proved that it can be expressed in terms of rational functions of $h^2$ and $\frac{1}{\sqrt{1-h^2}}$ as follows
\begin{eqnarray}
L_k(h)=\begin{cases}
P_{[\frac{-k}{2}]}(h^2),\qquad\qquad\quad\ \ k\leq0,\\
P_{[\frac{k-1}{2}]}(h^2)(1-h^2)^{\frac12-k},\ k\geq1,
\end{cases}
\end{eqnarray}
where $P_l$ are polynomials of exact degree $l$. In particular,
\begin{eqnarray}\begin{aligned}
L_0(h)=L_{-1}(h)=2\pi, L_{-2}(h)=\pi(2+h^2),L_1(h)=2\pi(1-h^2)^{-\frac12},L_2(h)=2\pi(1-h^2)^{-\frac32}.
\end{aligned}\end{eqnarray}
Now, we can proceed to calculate the derivatives of integrals $A_k(h)$. Some routine manipulation using (2.9) yields
\begin{eqnarray}\begin{aligned}
&\frac{d A_0(h) }{d h}=\frac{d A_1(h) }{d h}=\frac{2\pi}{h}-\frac{2\pi}{h\sqrt{1-h^2}},\\&\frac{d A_2(h) }{d h}=\frac{\pi}{h}-\frac{2\pi}{h^3}-\frac{2\pi}{h^3\sqrt{1-h^2}}.
\end{aligned}\end{eqnarray}

Next, we study $I_{i,j}(h)$ for $N\geq3$. It is apparent from Binomial Theorem that
\begin{eqnarray}\begin{aligned}
\int_0^{2\pi}\frac{\sin^k\theta}{(1-h\sin\theta)^n}d\theta=&h^{-k}\int_0^{2\pi}\frac{(1-(1-h\sin\theta))^k}{(1-h\sin\theta)^n}d\theta\\=&h^{-k}\sum\limits_{i=0}^k(-1)^i\binom{k}{i}L_{n-i}(h).
\end{aligned}\end{eqnarray}
It follows from (2.11) that
\begin{eqnarray*}\begin{aligned}
I_{i,j}(h)=&\sum\limits_{k=0}^Nd_kh^{N-2}\int_{-\pi}^\pi\frac{\sin^k\theta}{(1-h\sin\theta)^{N-2}}d\theta\\
=&\sum\limits_{k=0}^Nd_kh^{N-k-2}\sum\limits_{i=0}^k(-1)^i\binom{k}{i}L_{N-2-i}(h),
\end{aligned}\end{eqnarray*}
which leads to
\begin{eqnarray}\begin{aligned}
J_N(h)&=h^{-2}\sum\limits_{j=0}^N\varphi_j(h)L_{j-2}(h)\\&=h^{-2}\big[(1-h^2)^{\frac52- N}\varphi_{2N-3}(h)+\varphi_2(h)\big]
,\ N\geq3,
\end{aligned}\end{eqnarray}
in view of (2.8), where $\varphi_j(h)$ are polynomials of degree $j$.

In conclusion, on account of (2.4), (2.6), (2.10) and (2.12), one derives
\begin{eqnarray*}\begin{aligned}
\frac{dI_1(h)}{dh}=&\frac{d}{dh}\Big(\sum\limits_{N=1}^{[\frac{n+1}{2}]}J_N(h)-C_\delta\Big)\\
=&\frac{d J_1(h)}{dh}+\frac{d J_2(h)}{dh}+\frac{d}{dh}\Big(\sum\limits_{N=3}^{[\frac{n+1}{2}]}J_N(h)\Big)\\
=&h^{-3}(1-h^2)^{\frac32-[\frac{n+1}{2}]}\big[\phi_{2[\frac{n+1}{2}]-1}(h)+\psi_{2[\frac{n+1}{2}]-1}(h)\sqrt{1-h^2}\big],
\end{aligned}\end{eqnarray*}
where $\phi_j(h)$ and $\psi_j(h)$ are polynomials of degree $j$.
As a result, this together with the Rolle's Theorem ,one can obtain that  $I(h)$ has at most $4[\frac{n+1}{2}]+1$ zeros on (0,1) owing to (2.4) and the above equality.

\section{Discussion}
 \setcounter{equation}{0}

This study aims to get the number of limit cycles for the cubic isochronous differential system $S_2^*$ under the perturbations of $n$th-degree polynomials. The key to addressing this problem lies in our use of the first integral of (1.4) to transform the Abelian integral $I(h)$ into the form (2.1).
This allows us to apply polar coordinate transformations and, in conjunction with existing methods, estimate the upper and lower bounds of the number of limit cycles.
Thus, this paper, along with reference \cite{GLLZ}, the upper bounds of number of limit cycles of $S_1^*$ and $S_2^*$ under polynomial smooth perturbations have been obtained. But the bifurcation of limit cycles of $S_3^*$ and $S_4^*$ are still open. This is because their first integrals are quite involved. In the future, these challenges may be overcome by developing more powerful and innovative methods, which is precisely the problem we are currently addressing.

 \vskip 0.2 true cm

\noindent
{\bf CRediT authorship contribution statement}
 \vskip 0.2 true cm

\noindent
Jihua Yang: Conceptualization,  Funding acquisition, Investigation, Methodology, Resources, Supervision, Writing-original draft,	Writing-review $\&$ editing.
Qipeng Zhang: Validation, Writing-original draft.
\vskip 0.2 true cm

\noindent
{\bf Conflict of Interest}
 \vskip 0.2 true cm

\noindent
The authors declare that they have no known competing financial interests or personal relationships that could have
appeared to influence the work reported in this paper.
\vskip 0.2 true cm

\noindent
{\bf Data Availability Statement}
 \vskip 0.2 true cm

\noindent
No data was used for the research in this article. It is pure mathematics.
\vskip 0.2 true cm

\noindent
{\bf Acknowledgment}
 \vskip 0.2 true cm

\noindent
Supported by the National Natural Science Foundation of China(12161069) and the Ningxia Natural Science Foundation of China(2022AAC05044).

\end{document}